\documentclass[10pt]{article}

\usepackage[latin1]{inputenc} 
\usepackage[T1]{fontenc}  

\usepackage{graphicx}
\usepackage{verbatim}
\usepackage{amssymb}
\usepackage{amsmath}
\usepackage{amsthm} 
\usepackage[all]{xy}
\usepackage{amsfonts}

\usepackage[francais]{babel}

\newenvironment{disarray}%
 {\everymath{\displaystyle\everymath{}}\array}%
 {\endarray}

\newtheorem{theo}{\textbf{Théorème}}[section]
\newtheorem{lemme}[theo]{\textbf{Lemme}}
\newtheorem{pte}[theo]{\textbf{Proposition}}

\newtheorem{corollaire}[theo]{\textbf{Corollaire}}
\newtheorem{conjecture}[theo]{\textbf{Conjecture}}
\newtheorem{déf}[theo]{\textbf{Définition}}

\newcommand{\X}{\boldsymbol{X}}
\newcommand{\G}{\mathbb{G}}
\newcommand{\LL}{\mathbb{L}}
\newcommand{\FF}{\mathbb{F}}
\newcommand{\Q}{\mathbb{Q}}
\newcommand{\QQ}{\overline{\mathbb{Q}}}

\newcommand{\K}{\mathbb{K}}
\newcommand{\Z}{\mathbb{Z}}

\newcommand{\N}{\mathbb{N}}
\newcommand{\E}{\mathbb{E}}
\newcommand{\F}{\mathbb{F}}

\newcommand{\M}{\mathcal{M}}
\newcommand{\A}{\mathcal{A}}
\newcommand{\HH}{\mathcal{H}}
\newcommand{\p}{\mathcal{P}}
\newcommand{\OL}{\mathcal{O}_{\LL}}

\newcommand{\OK}{\mathcal{O}_{\mathbb{K}}}

\newcommand{\OO}{\mathcal{O}}

\newcommand{\W}{W}

\DeclareMathOperator{\codim}{codim}
\DeclareMathOperator{\Gal}{Gal}

\newcommand{\oo}{\Omega}
\newcommand{\zz}{\boldsymbol{\zeta}}

\newcommand{\xx}{\boldsymbol{x}}
\newcommand{\yy}{\boldsymbol{y}}
\newcommand{\al}{\boldsymbol{\alpha}}

\newcommand{\w}{\omega}
\newcommand{\g}{\gamma}
\newcommand{\z}{\zeta}

\newcommand{\D}{\Delta(W)}
\newcommand{\mm}{\hat{\mathcal{\mu}}_{ess}}
\newcommand{\eps}{\varepsilon}
\newcommand{\ro}{\rho}
\newcommand{\dd}{\delta}

\numberwithin{equation}{section}

\begin{document}

\title{\textbf{Problème de Lehmer relatif dans un tore :
 cas des hypersurfaces}}
\author{\sc{Emmanuel DELSINNE}\footnote{{\it UMR 6139 (CNRS), Laboratoire
    de Mathématiques Nicolas Oresme,
Université de Caen, 
BP 5186, 14032 Caen Cedex }
 (delsinne@math.unicaen.fr)}}

\date{}

\maketitle

\renewcommand{\abstractname}{Abstract}

\begin{abstract}
We tackle the ``relative'' Lehmer problem on algebraic
  subvarieties of a multiplicative torus. Generalizing a theorem of F.~Amoroso and
  U.~Zannier, we give a lower bound for the normalized height of a non
  torsion hypersurface in terms of its obstruction index over
  $\Q^{ab}$, the maximal abelian extension of $\Q$. We prove up to
  $\eps$ the sharpest conjecture that can be formulated.
\end{abstract}

\renewcommand{\abstractname}{Résumé}

\begin{abstract}
  Nous abordons le problème de Lehmer <<~relatif~>> pour les
  sous-variétés algébriques d'un tore multiplicatif. Généralisant
  un théorème de F.~Amoroso et U.~Zannier, nous montrons que la
  hauteur normalisée d'une hypersurface qui n'est pas de torsion est
  minorée en fonction de son indice d'obstruction sur $\Q^{ab}$,
  l'extension abélienne maximale de $\Q$. La minoration ainsi obtenue
  correspond à un $\eps$-près à la conjecture la plus précise que l'on
  peut formuler dans ce cadre.

\end{abstract}

\section{Introduction}

Nous nous proposons ici de poursuivre l'étude des minorations de la
hauteur normalisée des sous-variétés d'un tore amorcée par {\sc
  F.~Amoroso} 
et {\sc S.~David} dans \cite{Am-Da99}, \cite{Am-Da00}, \cite{Am-Da03} et
\cite{Am-Da04}. Soit $n$ un entier naturel non nul.
Nous considérons le plongement ``naturel'' de $\G_m^n$ dans
$\mathbb{P}_n$. La hauteur (normalisée) d'un point $\al =
(\alpha_1,\dots, \alpha_n) \in \G_m^n $ est donc la hauteur de Weil
logarithmique et absolue (avec la norme du sup aux places
archimédiennes) $\hat{h}(\al)$ du point projectif
$(1:\alpha_1:\dots:\alpha_n)$.
{\sc P.~Philippon} (\cite{Ph1},\cite{Ph2},\cite{Ph3})  d\'efinit
la hauteur normalis\'ee d'une sous-vari\'et\'e $V$ de $\G_m^n$ par~:
\begin{equation*}
\hat{h}(V)=\lim_{m\rightarrow+\infty}
\frac{h([m]V)\deg(V)}{m\deg([m]V)}\enspace,
\end{equation*}
o\`u $h(V)$ (respectivement $\deg(V)$) est la hauteur projective
(respectivement le degr\'e) de
l'adh\'erence de Zariski de $V$ dans $\mathbb{P}_n$. {\sc L.~Szpiro} a
\'egalement introduit
le {\it minimum essentiel} de $V$, not\'e $\mm(V)$, comme la borne
inf\'erieure des nombres
r\'eels $\theta>0$ tels que l'ensemble des points $P\in V({\QQ})$ de
hauteur normalis\'ee
born\'ee par $\theta$ soit Zariski-dense dans $V$. Si $V$ est
$\QQ$-irr\'eductible, on dispose alors de la
relation suivante, montr\'ee dans~\cite{zhang1995a} et \cite{zhang1995b}~:
\begin{equation*}
\frac{\hat{h}(V)}{(\dim(V)+1)\deg(V)}
\leq\mm(V)\leq\frac{\hat{h}(V)}{\deg(V)}\enspace,
\end{equation*}
Le minimum essentiel et la hauteur normalis\'ee ont la propri\'et\'e
remarquable suivante, montr\'e
encore par {\sc S.~Zhang} ({\em confer} \cite{zhang1992})~:
$\mm(V)=0$ (et donc $\hat{h}(V)=0$)
si et seulement si $V$ est une vari\'et\'e de torsion ({\it i.e.} une
r\'eunion de translat\'es de
sous-tores de $\G_m^n$ par des points de torsion).

Il est donc naturel de chercher \`a minorer le minimum essentiel (ou
la hauteur normalis\'ee) d'une
vari\'et\'e qui n'est pas de torsion (ou bien imposer des conditions
g\'eom\'etriques portant sur la dimension du
stabilisateur de $V$).




  Une telle minoration va dépendre
 des caractéristiques géométriques de la variété, par exemple son
 degré. Cependant, si l'on n'impose aucune condition géométrique sur la
 variété, il faudra également tenir compte de son corps de définition.
 En effet, soit $H$ un sous-groupe de $\G_m^n$ et soit $\al_i$ une
 suite de points de non-torsion dont la hauteur tend vers $0$ (par
 exemple $\al_i = (2^{1/i},\dots,2^{1/i})$). Alors les variétés
 $V_i:=H\al_i$ ont toutes même degré $\deg(H)$ mais leur minimum
 essentiel $\mm(V_i) \leq h(\al_i)$ converge vers $0$.

  Le problème consistant à trouver les meilleures bornes inférieures
  pour le minimum essentiel des sous-variétés de $\G^m_n$ est une
  généralisation d'une célèbre question de {\sc D.~H.~ Lehmer} :
  existe-t-il un constante $c>0$ telle que pour tout 
  nombre algébrique $\alpha$ de degré $d$ qui n'est pas une racine de
  l'unité on ait $h(\alpha) \geq c/d $ ?
  Si l'on ne suppose rien de plus sur $\alpha$, c'est la meilleure
  minoration possible, étant donné que $h(2^{1/d}) = (\log 2) / d $. 
  Dans cette direction, le meilleur résultat à ce jour est un résultat de
  {\sc E.~Dobrowolski} :

  \begin{theo}
   Il existe une constante $c>0$ tel que pour tout nombre algébrique
   $\alpha$ de degré $d(\geq 2)$ qui n'est pas une racine de l'unité :
  $$
  h(\alpha) \geq \frac{c}{d} \left ( \frac{\log \log d }{\log d} \right )^3.
  $$ 
  \end{theo}

Cependant {\sc F. Amoroso} et {\sc U. Zannier} ont montré dans
\cite{Am-Za} que l'on a le même type de minoration en remplaçant le
degré total de $\alpha$ ({\it i.e.} $d = [\Q(\alpha):\Q]$) par 
le degré <<non abélien>> de 
$\alpha$ ({\it i.e.} $[\Q^{ab}(\alpha):\Q^{ab}]$ où $\Q^{ab}$ désigne
l'extension abélienne maximale de $\Q$). Notre but est
de généraliser ce résultat en dimension supérieure.

 Dans les problèmes de minoration en dimension supérieure, l'invariant
 le plus fin qui puisse tenir compte de la nature <<arithmétique>>
 d'une variété est 
 l'\emph{indice d'obstruction}.
 Quelques notations sont nécessaires
 avant d'introduire cet invariant.
Nous fixons donc $\QQ$ une clotûre algébrique de $\Q$ et nous noterons
$\G_m^n$ pour $\G_m^n (\QQ)$.
Soit $V$ une sous-variété de $\G_m^n$ et soit $\K$ un sous-corps de
$\QQ$. Nous utiliserons les notations suivantes : $\Q(V)$ désignera le
corps de définition de $V$, $\K(V)$ le corps $\K \cdot \Q(V)$ et
$\overline{V}^{\K}$ la variété définie par :
$$
\overline{V}^{\K} : = \bigcup_{\sigma \in \Gal (\QQ / \K) } \sigma V.
$$
Remarquons que $\deg \overline{V}^{\K} = [\K(V) : \K] \deg V$.

\begin{déf}
  On appelle \emph{indice d'obstruction de $V$ sur $K$} l'entier 
  $\omega_{\K} (V) $ défini par :
  $$
\omega_{\K} (V) = \min_{Z \supseteq V  \atop \codim Z= 1} \{ \deg
(\overline{Z}^{\K} )   \}.  
$$ 
\end{déf}
\noindent Par exemple, si $V$ est une hypersurface de $\G_m^n$, $\w_\K(V) =
[\K(V) : \K] \deg V.  $

{\sc F. Amoroso} et {\sc S. David} énoncent alors une conjecture
généralisant le problème de Lehmer en dimension supérieure et
obtiennent dans cette direction le résultat suivant,
analogue du théorème de {\sc Dobrowolski} en dimension supérieure :

\begin{theo}
Soit $n$ un entier naturel non nul.
Soit $\W$ une sous-variété géométriquement irréductible de $\G^n_m$ de
codimension $k$
qui n'est contenue dans aucune sous-variété de torsion.
Alors
$$
\mm(\W) \geq \frac{c(n)}{{\w}_{\Q} (\W)}  \log 
( 3 \w_{\Q} (\W) ) ^{-\lambda(k)}  
$$
où $c(n)$ et $\lambda(k)$ sont des constantes strictement positives ne
dépendant respectivement que de $n$ et $k$.
\end{theo}

Nous pouvons ainsi énoncer la conjecture <<abélienne>> analogue :

\begin{conjecture}
Soit $n$ un entier naturel non nul.
Soit $\W$ une sous-variété géométriquement irréductible de $\G^n_m$
qui n'est contenue dans aucune sous-variété de torsion.
Soit $\LL$ une extension abélienne de $\Q$.
Alors $$\mm(\W) \geq \frac{c(n)}{{\w}_{\LL} (\W)} \  $$
où $c(n)$ est une constante strictement positive ne dépendant que de
$n$.  
\end{conjecture}

Dans cette direction,  en <<combinant>> les techniques de
\cite{Am-Da00} et \cite{Am-Za}, nous obtenons le résultat suivant
concernant les variétés de codimension 1:

\begin{theo}\label{th}
Soit $n$ un entier naturel non nul.
Soit $\W$ une hypersurface géométriquement irréductible de $\G^n_m$
qui n'est pas de torsion. Soit $\LL$ une extension abélienne de $\Q$.
Alors $$\mm(\W) \geq \frac{c(n)}{{\w}_{\LL} (\W)} \bigg ( \frac{\log 
 2 \w_{\LL} (\W)}{\log \log 5 \w_{\LL} (\W)}\bigg )^{-(1+6(n+1))} $$
où $c(n)$ est une constante strictement positive ne dépendant que de $n$.
\end{theo}

Remarquons également que pour les hypersurfaces, on peut attaquer ce
problème d'un point de vue différent. En effet, dans \cite{Am-Da03},
{\sc F.~Amoroso} et {\sc S.~David} obtiennent, en introduisant des
hypothèses supplémentaires, une minoration uniquement de <<nature
géométrique>> . 

\begin{theo}
Soit $n$ un entier naturel non nul.
Soit $\W$ une sous-variété géométriquement irréductible de $\G^n_m$ de
codimension $k$
qui n'est contenue dans aucun translaté de sous-tore.
Alors
$$
\mm(\W) \geq \frac{c(n)}{{\w}_{\QQ} (\W)}  \log 
( 3 \w_{\QQ} (\W) ) ^{-\lambda(k)}  
$$
où $c(n)$ et $\lambda(k)$ sont des constantes strictement positives ne
dépendant respectivement que de $n$ et $k$.
\end{theo}

Ce résultat, appliqué aux hypersurfaces, nous indique que si $W$ est
une hypersurface géométriquement irréductible qui n'est pas le translaté
d'un sous-tore alors :
\begin{equation}\label{geom}
 \mm(\W) \geq \frac{c(n)}{\deg \W}  \log 
( 3 \deg \W ) ^{-81}  
\end{equation}
Pour obtenir un résultat du type Théorème \ref{th}, il suffit donc de
ne s'intéresser qu'aux hypersurfaces qui sont des translatés de
sous-tores par des points d'ordre infini. Mais dans ce cas on peut
facilement se ramener au cas du théorème principal de \cite{Am-Za} et
obtenir :

\begin{theo}\label{sous-tore}
Soit $n$ un entier naturel non nul.
Soit $\W$ une hypersurface géométriquement irréductible de $\G^n_m$
qui est le translaté d'un sous-tore par un point d'ordre infini.
Alors
$$
\mm(\W) \geq \frac{c}{n \cdot \w_{\LL}(W)} \bigg ( \frac{\log 
( 2 [\LL(W):\LL] )}{\log \log (5 [\LL(W):\LL] )} \bigg )  ^{-13}  
$$
où $c$ est une constante strictement positive.
\end{theo}

{\it Démonstration. }
Soient $H$ un sous-tore géométriquement irréductible de codimension 1
et $\al$ un point d'ordre infini tels que $W=\al H$. En tant que
sous-tore de $\G_m^n$ de codimension 1, $H$ est donné par
une équation du type $\X^{\boldsymbol{\lambda}} = 1 $ avec
${\boldsymbol{\lambda}} \in \Z^n$. Si on pose $\boldsymbol{\mu} =
(\mu_1, \dots, \mu_n) \in \N^n$ avec $\mu_i = \max ( 0, -\lambda_i) $ et $
\alpha = \al^ {\boldsymbol{\lambda}} \in \QQ$, alors le polynôme $F(\X) = \X^
{\boldsymbol{\mu }} ( \X^{\boldsymbol{\lambda}} -\alpha)$ est une
équation de $W$ et il est clair que $\LL(W) = \LL(\alpha)$. D'une
part, comme W est une hypersurface, on a $\hat h (W) = \hat h (F)
$. D'autre part, un simple changement de variables dans les calculs de
la mesure de Mahler de $F$ nous donne $\hat h (F) = \hat h (X-\alpha)
= \hat h (\alpha)$ (pour les liens entre la hauteur normalisée des
hypersurfaces et la mesure de Mahler de leurs équations, voir par
exemple, \cite{Ph1}, section 2, partie B). Ainsi en appliquant le
théorème de \cite{Am-Za}, on obtient :
$$
\hat h (W) \geq \frac{c}{  [\LL(W):\LL] } \bigg ( \frac{\log 
( 2 [\LL(W):\LL] )}{\log \log (5 [\LL(W):\LL] )} \bigg )  ^{-13}  
$$
On conclut alors en utilisant l'inégalité de Zhang. $\hfill \square$

\bigskip
En utilisant l'inégalité (\ref{geom}) et le théorème \ref{sous-tore}
on peut donc obtenir un résultat du type théorème \ref{th} avec une
constante absolue comme exposant du terme en <<$\log$>>. Cependant, cette
approche du problème ne fonctionne plus en codimension supérieure et il faudra
sans doute privilégier un raisonnement combinant les techniques de
\cite{Am-Za} et \cite{Am-Da99}. L'objet de ce papier est donc de
démontrer le théorème \ref{th} par cette voie.


\section{Notations et résultats préliminaires}

Soit $n$ un entier naturel non nul.
Soient $\xx, \ \yy \in \G_m^n$ et soit $m \in \N^*$.
On notera :
$$ \xx  \yy = ( x_1y_1, \cdots, x_ny_n) \ \text{ et } \
[m]\xx=(x_1^m, \cdots , x_n^m).$$
On désignera par $\ker[m]$ le noyau du morphisme de multiplication par
$m$ dans $\G_m^n$, {\it i.e.} l'ensemble des points dont les
coordonnées sont des racines $m$-ièmes de l'unité. Si $V$ est une
sous-variété de $\G_m^n$, on notera $G_V$ son stabilisateur :
$$G_V = \{ \xx \in \G_m^n , \ \xx  V = V \} = \bigcap_{ \yy \in V }
\yy^{-1}   V $$
et $G_V^0$ la composante neutre de $G_V$. Le stabilisateur de $V$
possède les propriétés suivantes : 
$$
\dim (G_V) \leq \dim (V) \ \text{ et } \ \deg(G_V) \leq
{\deg(V)}^{\dim(V) + 1}.
$$
Par ailleurs, si $W$ est une sous-variété stricte et
\emph{géométriquement irréductible} de $\G_m^n$, le degré de son image
par le morphisme de multiplication par $m$ vérifie :
$$
\deg([m]W) = \frac{m^{\dim(W)} \deg (W) }{|\ker[m] \cap G_W|}
= \frac{m^{\dim(W)-\dim{G_W}} \deg (W) }{|\ker[m] \cap (G_W/G_W^0)|}.
$$
où l'on a encore noté $\ker [m]$ le noyau de la multiplication par $m$
dans $\G^n_m/G_W^0$.
On pourra trouver une démonstration de ces résultats dans
\cite{Am-Da99} et \cite{Hindry}. Enfin, nous aurons besoin du lemme
suivant :

\begin{lemme}\label{zeta}
Soit $W$ une sous-variété de $\G_m^n$ géométriquement irréductible.
Soient $\K$ un corps de nombres, $p$ un nombre premier et $\z_p$
une racine primitive $p$-ième de l'unité. Alors l'extension 
$$\K([p]W,\z_p) \subseteq \K(W,\z_p) $$
est abélienne de degré une puissance de $p$. De plus, si 
$\K(W,\z_p) = \K([p]W,\z_p)$, 
il existe $\zz \in \ker[p]$ tel que $\K(\zz W) = \K([p]W)$.
\end{lemme}

\bigskip

{\it Démonstration. }
Soit $\tau \in \Gal(\QQ / \K([p]W, \z_p))$. Montrons que
$\K(W,\z_p)$ est globalement stable sous-l'action de $\tau$. Pour
cela, il suffit de montrer que $\tau (W) $ est définie sur $
\K(W,\z_p) $. On a :
$$
[p] \tau (W) = \tau([p] W ) = [p]W .
$$
Il existe donc $\boldsymbol{\xi} \in \ker[p] $ tel que $
\tau (W) = \boldsymbol{\xi} W$ et $\tau (W) $ est définie sur $
\K(W,\z_p) $. L'extension 
$\K([p]W,\z_p) \subseteq \K(W,\z_p) $ est donc galoisienne. D'autre
part, si l'on considère l'application :
$$
\begin{array}{rrcl}
\phi: & \Gal ( \K(W,\z_p) / \K([p]W,\z_p) )  & \longrightarrow &   
    \left( ^{ \ker [p]}/_{ \ker [p] \cap G_{W}} \right ) \\[4mm]
 & \tau & \longmapsto &  \bar{ \boldsymbol{\xi}}
\end{array}
$$
on vérifie aisément que $\phi$ est bien définie et que c'est un
morphisme injectif. Ainsi $\Gal ( \K(W,\z_p) /
\K([p]W,\z_p) )  $ est isomorphe à son image par $\phi$ et donc est
abélien. La première partie du lemme et donc démontrée, passons à la
seconde.

Remarquons d'abord que, par hypothèse, 
$$
\K([p] W) \subseteq \K(W) \subseteq \K(W,\z_p) = \K([p]W, \z_p).$$
Si $\K(W) = \K([p]W)$ le résultat est trivial. Supposons donc que
$\K([p] W) \subsetneq \K(W) $. Soit $\sigma$ un générateur du groupe
cyclique 
$$
G = \Gal ( \K([p]W,\z_p) / \K([p]W) )
$$
et notons $\tilde \sigma$ un de ses prolongements à $\overline{\Q}$.
Comme $\tilde \sigma $ fixe $\K([p] W )$ et {\it a fortiori}
$\Q([p] W )$, on a~:
$$
 [p] \tilde \sigma (W) = \tilde \sigma ([p] W) = [p] W.
$$
Il existe donc $\boldsymbol{\xi} \in \ker[p] $ tel que $\tilde
\sigma (W) = \boldsymbol{\xi} W$.

Montrons que $\sigma(\boldsymbol{\xi} ) \neq \boldsymbol{\xi}$.
Si $\boldsymbol{\xi} = (1,\dots,1) $ alors $\tilde \sigma (W) = W$ et
$\Q(W)$ est stable sous l'action de $G$; on
en déduit que $\K(W) = \K([p]W) $. Par ailleurs, si $\boldsymbol{\xi}
\neq (1,\dots,1)$ et $\sigma \boldsymbol{\xi} = \boldsymbol{\xi}$,
alors $\K(\boldsymbol{\xi}) = \K(\z_p)$ est stable sous l'action de
$G$~;~il s'en suit que $G$ est réduit à l'identité et $\K([p] W
)=\K([p]W,\z_p)$~;~{\it a fortiori} on a encore $ \K(W) = \K([p]W)
$. Dans les deux cas, on obtient une contradiction avec
l'hypothèse $\K(W) \neq \K([p]W) $.

On a donc $\sigma (\boldsymbol{\xi}) = \boldsymbol{\xi}^{\lambda} $,
avec $\lambda \in \Z$ et $\lambda \not \equiv 1  \mod p$. Soit $u$ une
solution de la congruence 
$$
(\lambda - 1)u + 1 \equiv 0  \mod p 
$$
et soit $\zz := \boldsymbol{\xi}^u$. On a~:
$$
\tilde \sigma ( \zz W ) = \sigma (\zz ) \tilde \sigma (W) =
\boldsymbol{\xi}^{\lambda u +1} W = \zz W ,
$$
ce qui montre que $\Q(\zz W)$ (donc $\K(\zz W)$) est stable sous
l'action de $G$, et ainsi $\K(\zz W ) \subseteq \K([p]W)$. D'autre
part, $[p](\zz W) = [p] W $, donc ces deux corps sont égaux, ce qui
achève la démonstration. $\hfill \square$

\section{Réductions}

Soit $\LL$ une extension abélienne de $\Q$.  D'après le
théorème de Kronecker-Weber, $\LL$ est contenu dans une extension
cyclotomique de $\Q$. Soit $m \in {\N}^* $ minimal tel que $\LL
\subseteq  
\Q(\z_m)$. Si $p$ est un nombre premier,
on note $e_p(\LL)$ son indice de ramification dans $\LL$ et
$\tilde{e_p} (\LL)$ la puissance maximale de $p$ divisant $m$. 
On définit également
$$\tilde {e} (\LL) = \sum_{p \text{ premier }} (\tilde{e_p} (\LL) -
1). $$
Remarquons que si $\LL' \subseteq \LL$ sont deux extensions abéliennes
de $\Q$ alors  $\tilde {e} (\LL') \leq \tilde {e}
(\LL)$.

Pour démontrer le théorème \ref{th}, on raisonne par l'absurde. 
Supposons donc qu'il existe $\LL$ une extension abélienne de $\Q$ et
$\W$ une hypersurface  géométriquement irréductible de $\G_m^n$ non de
torsion tels que le théorème \ref{th} soit faux :
\begin{equation}\label{1}
\mm(\W) < \frac{c(n)}{{\w}_{\LL} (\W)} \bigg ( \frac{\log
 2  \w_{\LL} (\W)}{\log \log 5 \w_{\LL} (\W)}\bigg )^{-(1+6(n+1))}.
\end{equation}

Nous pouvons supposer de plus que le degré $\dd : = [\LL(\W): \LL]$
  est minimal
  dans (\ref{1}), {\it i.e.} pour toute hypersurface géométriquement
  irréductible  $\W'$ qui n'est pas
  de torsion et telle qu'il existe une extension abélienne $\LL'$ de 
$\Q$ vérifiant $[\LL'(\W'): \LL'] <  \dd $, on a :
\begin{equation}\label{2}
\mm(\W') \geq \frac{c(n)}{{\w}_{\LL'} (\W')} \bigg ( \frac{\log 2
  \w_{\LL'} (\W')}{\log \log 5 \w_{\LL'} (\W')}\bigg )^{-(1+6(n+1))}.
\end{equation}
Remarquons ensuite que la fonction
$$
t \mapsto t \cdot \left (\frac{\log (2 (\deg \W)  t) } {
\log \log (5 (\deg \W) t) } \right )^{1+6(n+1)} 
$$
est croissante sur $[1; +\infty[$. De plus, pour tout $ \zz \in
(\G^n_m)_{\text{tors }}$, on a $\deg (\zz W) = \deg(W)$ et
$\mm(\zz\W) = \mm(\W)$. En particulier, (\ref{1}) et (\ref{2})
impliquent que 
pour tout  $\zz \in (\G^n_m)_{\text{tors}}$ et toute extension
abélienne $ \LL'$  de $ \Q$, on a :  
\begin{equation}\label{3}
\quad[\LL'(\zz \W): \LL'] \geq \dd.
\end{equation}

Enfin, soit $\A$ l'ensemble des extensions abéliennes $\LL'$ de $\Q$
telles qu'il existe $\zz \in (\G^n_m)_{\text{tors}}$ vérifiant
$[\LL'(\zz \W): \LL'] = \dd  $. Soit
 $$ \tilde{e} : =
\min_{\LL \in  \A} \tilde{e}
(\LL). $$
Quitte à remplacer $\W$ par $\zz \W$ pour un certain $\zz \in
(\G^n_m)_{\text{tors}}$ et $\LL$ par $\LL' \in \A$  nous pouvons 
supposer que $\LL$ vérifie les deux conditions suivantes :

\begin{equation}\label{4}
 [{\LL} (\W): {\LL}] = \dd
\end{equation}
\begin{equation}\label{5}
\tilde{e}({\LL}) = \tilde {e}
\end{equation}
De plus, par un argument galoisien, nous avons
le diagramme suivant~:

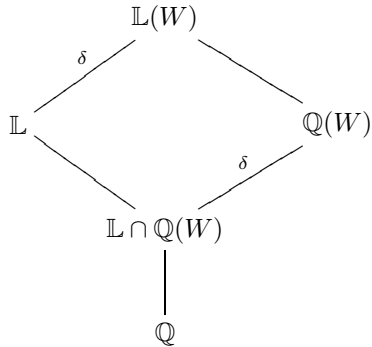
  $$
  \xymatrix{
               & \LL(W) \ar@{-}[ld]_{\dd} \ar@{-}[rd]&    \\ 
\LL \ar@{-}[rd]&                               &\Q(W) \ar@{-}[ld]_\dd  \\
               & \LL \cap \Q(W)                & \\
               & \Q \ar@{-}[u]                 &
  }
  $$
\'Etant donné que $\LL \cap \Q(W)$ est une extension abélienne de $\Q$ et que
$$\tilde{e}( \LL\cap\Q(W)) \leq \tilde{e} (\LL) = \tilde{e},$$
on a $\tilde{e}( \LL \cap \Q(W)) = \tilde{e}$.
Cela nous permet de supposer, quitte à remplacer $\LL$ par $\LL \cap
\Q(W) $ que $\LL \subseteq \Q(W)$, {\it i.e}~:
\begin{equation}\label{6}
\Q(\W)=\LL(\W).
\end{equation}
Remarquons enfin que
l'on peut également supposer que~: 
\begin{equation}\label{7}
\forall \zz \in (\G^n_m)_{\text{tors}},\quad  \Q(\zz \W ) \subseteq \Q
(\W)
\Rightarrow \Q(\zz \W ) =  \Q (\W).
\end{equation}
En effet, s'il existe $\zz$ tel que  $\Q(\zz W) \subsetneq \Q(W)$,
nous avons $\LL(\zz W) \subseteq \LL(W)$, ce qui implique
nécessairement (par \ref{3}) $ \LL(\zz W) = \LL(W) = \Q(W) $.
Nous avons ainsi le diagramme suivant :

 $$
  \xymatrix{
               & \Q(W) \ar@{-}[ld]_{\dd} \ar@{-}[rd]&    \\ 
\LL \ar@{-}[rd]&                               &\Q(\zz W) \ar@{-}[ld]_\dd  \\
               & \LL \cap \Q(\zz W)            & \\
               & \Q \ar@{-}[u]                 &
  }
  $$
Par le même argument, nous pouvons donc remplacer $W$ par $\zz W$ et
$\LL$ par $\LL \cap \Q(\zz W)$. Nous pouvons itérer ce procédé,
jusqu'à obtenir (\ref{7}) (nombre d'itérations fini car le degré décroit
strictement à chaque étape).

Ainsi, nous considérons désormais une hypersurface géométriquement
irréductible $\W$ qui n'est pas de torsion et une extension abélienne
$\LL$ de $\Q$ contenue dans $\Q(\W)$ qui satisfont (\ref{1}),
(\ref{2}), (\ref{3}), (\ref{4}), (\ref{5}), (\ref{6}) et (\ref{7}).


\bigskip

\noindent {\bf Notations.}
Soit $m \in \N^*$. Dans toute la suite, nous noterons par
$V_m$ la variété définie par :
\begin{equation}
V_m = \bigcup_{\sigma \in \Gal (\QQ / \LL) } [m] \sigma (W) = [m]
\overline{W}^{\LL} .
\end{equation}
On définit également un ensemble de premiers ``exceptionnels'' :
$$ E_{\text{exc}}(V_1) = \{  p \text{ premiers} , p \mid
|G_{V_1}/G_{V_1}^0| \},$$
et on note $\p$ son complémentaire dans l'ensemble des nombres
premiers.
Enfin, nous noterons $s$ la dimension du stabilisateur de $W$.

\begin{lemme}\label{i}
Soit $p \in \p $.
Alors :
$$
\begin{array}{rl}
i) & p \text{ ne divise pas } \ | G_\W/G^0_\W |, \\[2mm]
ii) &  \LL ([p]\W)=\LL ( \W), \\[2mm]
iii) & \deg( V_p ) = p^{n-1-s}\w_{\LL}(W). 
\end{array}
$$
\end{lemme}

\bigskip

{\it Démonstration. }$i)$ On montre facilement que $G_W
\subseteq G_{V_1}$ et $G_W^0 = G_{V_1}^0$. Ainsi $ G_\W/G^0_\W $ est un
sous-groupe de $G_{V_1}/G^0_{V_1}$ et $p$ ne divise pas $|G_\W/G^0_\W
|$. 

$ii)$ Il suffit de montrer que $\LL([p]\W,\z_p) = \LL(\W,\z_p)$ : le
lemme \ref{zeta} nous indique alors l'existence de $\zz \in
\text{ker}([p])$ tel que $\LL(\zz\W)=\LL([p]\W)$. Les conditions
(\ref{3}) et (\ref{4}) faite sur $\W$ nous donnent ainsi :
$$
\dd = [\LL(W):\LL] \geq [\LL([p]W):\LL] = [\LL(\zz\W) : \LL] \geq \dd.
$$
On en déduit que $\LL ([p]\W)=\LL ( \W) $.

Considérons l'extension abélienne $$\LL([p]\W,\z_p) \subseteq \LL(\W,\z_p)$$
et supposons qu'il existe un élément
$\sigma\neq \text{Id}$ dans $ \Gal( \LL(\W,\z_p) / \LL([p]\W,\z_p))  $ et
notons $\tilde \sigma$ un de ses prolongements à $\QQ$. On a 
$$
[p]\tilde \sigma W = \tilde \sigma [p] W = [p] W
$$
donc il existe $\boldsymbol{\xi} \in \ker [p] $ différent de
$(1,\dots,1)$ tel que $\tilde \sigma (W) = \boldsymbol{\xi} 
W$. Soit $\tau \in \Gal (\QQ / \LL)$; il existe $l \in \Z$ tel que que
$\tau^{-1} \boldsymbol{\xi} = \boldsymbol{\xi}^l$. On a alors, comme
$\sigma(\boldsymbol{\xi}) =
\boldsymbol{\xi}$ :
$$
\boldsymbol{\xi} \tau (W) = \tau (\boldsymbol{\xi} ^l W) = ( \tau
 \circ \tilde \sigma^l ) (W). 
$$
Donc $\boldsymbol{\xi} \in G_{V_1}$. Mais comme $G_{V_1}^0=G_W^0$ et $\sigma
\neq$ Id, on a $ \boldsymbol{\xi}
\not \in G_{V_1}^0$. On en déduit que $p$ divise
$|G_{V_1}/G^0_{V_1}|$, ce qui est absurde. On a donc bien
$\LL([p]\W,\z_p) = \LL(\W,\z_p)$ et le point $(ii)$ est établi.

$iii)$ Le point précédent nous assure que $[p]W$ et $W$ ont le même
nombre ($=[\LL(W):\LL]$) de conjugués au dessus de $\LL$. D'où :
$$
\deg(V_p) = [\LL(W) : \LL] \deg([p]W).
$$ 
Or, d'après $(i)$, $|\ker[p] \cap G_W/G_W^0| = 1$. Donc $\deg([p]W)=
p^{n-1-s} \deg W$, et la preuve du lemme est achevée.

 $\hfill \square$ 

\bigskip




\section{Lemmes pour l'extrapolation}

\begin{lemme}\label{Phi}
Soit $p$ un nombre premier. Soit $(p) = (\pi_1 \cdots
\pi_r)^{e_p(\LL)}$ la décomposition de $(p)$ dans $\OL$.
Alors il existe un élément $\Phi_p$
du groupe de Galois $\text{Gal}(\LL/\Q) $ tel que pour tout entier
algébrique $\g \in \LL$ :  
$$  \g^p - \Phi_p \g  \equiv 0 \ \mod \pi_1 \cdots
\pi_r  . $$
\end{lemme}

\bigskip

{\it Démonstration. }
Voir \cite{Am-Za}, Lemme 3.1. $\hfill \square$ 

\bigskip

\begin{lemme}\label{H_p}
Soit $p \in \p$. Alors il existe un sous-groupe $H_p$ de $\Gal (\LL /
\Q) $ d'ordre 
$$| H_p | \geq \min \{e_p(\LL), p \} $$
tel que pour tout entier
algébrique $\g \in \LL$, pour tout $\sigma \in H_p$, on ait~:  
\begin{equation}\label{mod p}
   \g^p - \sigma \g^p  \equiv  0 \  \mod  p \OL  .
\end{equation}
De plus, pour tout prolongement $\tau \in \Gal(\overline \Q / \Q)$ de
$\sigma  \in
H_p\backslash \{Id\}$, on a~:
 $$\tau [p]W \neq [p]  W. $$
\end{lemme}

\bigskip

{\it Démonstration. } Si $p$ n'est pas ramifié dans $\LL$ alors
$e_p(\LL) = 1$ et $H_p= \{ Id \}$ auquel cas le
lemme est trivial.
Si $p$ est ramifié dans $\LL$ alors $p$ est également ramifié dans
$\Q(\z_m)$, donc $p|m$. Notons $G_p : = \Gal ( \Q(\z_m) /
\Q(\z_{m/p}))$ qui est cyclique d'ordre $p$ si $p^2 | m$, d'ordre
$p-1$ sinon. Par minimalité de $m$, $\LL$ n'est pas stable sous
l'action de $G_p$ donc $G_p$ induit par restriction un sous-groupe
non trivial $H_p$ de $\Gal(\LL/ \Q)$. Si $p^2 | m $, alors
nécessairement $|H_p|=p$. Si $p^2 \nmid m$, alors $|H_p| \mid (p-1) $
et $|H_p| \geq e_p(\LL) $ car $p$ n'est pas ramifié dans
$\Q(\z_{m/p})$.

Soit $\g \in \LL$ un entier algébrique. En particulier, $\g$ est un
entier de $\Q(\z_m)$, donc s'écrit $\g = f(\z_m) $, avec $f
\in \Z[X]$. Soit $\sigma \in H_p$. Cet automorphisme est la
restriction à $\LL$
d'un certain $\tilde \sigma \in G_p$. Comme $\Q(\z_{m/p})$ est stable
par l'action de $\tilde \sigma$, on a $ \tilde \sigma ( \zeta_m^p ) =
\zeta_m^p  $.
On obtient ainsi, à l'aide du petit théorème de Fermat : 
$$\tilde  \sigma \g^p =\tilde \sigma f(\zeta_m)^p \equiv \tilde
\sigma f(\zeta_m^p)  = 
f(\zeta_m^p)  \equiv \g^p \quad \mod  p \Z[\z_m]. $$
Ce qui, par restriction à $\LL$, nous donne (\ref{mod p}).

Enfin, soient $\sigma \in H_p\backslash\{Id\}$ et $\tau \in
\Gal(\overline \Q / \Q)$ un prolongement de $\sigma$.
Supposons que $\tau [p]W  = [p] W $.
Ceci équivaut à dire que $\Q([p]W)$ est stable sous l'action de
$\tau$. 
Notons $\E$ le sous-corps de
$\LL$ fixé par $\sigma$.
On a alors $\Q([p] W) \cap \LL \subseteq \E$ donc $\E([p] W) \cap
\LL = \E$. Par un argument galoisien, les <<cotés>> opposés du
diagramme suivant ont même degré :

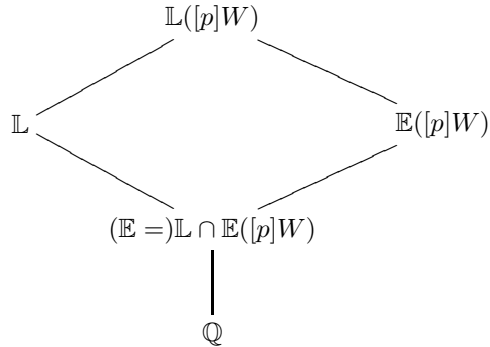
 $$
  \xymatrix{
               & \LL([p]W) \ar@{-}[ld] \ar@{-}[rd]&    \\ 
\LL \ar@{-}[rd]&                               &\E([p]W) \ar@{-}[ld]  \\
               & (\E =) \LL \cap \E([p]W)      &    \\
               & \Q \ar@{-}[u]                 &
  }
  $$
On en déduit :
$$
 [\LL([p]W) : \E([p]W) ] = [\LL:\E].
$$
De plus, par le lemme \ref{i} et l'égalité (\ref{6}), on a   
\begin{equation}\label{égalité}
 \LL([p]W)  =   \LL(W) = \Q(W) .
\end{equation}
Ainsi $[\Q(W) : \E([p]W) ] = [\LL:\E] $.
D'une part, comme $\E \subsetneq \LL$ (sinon $\sigma = Id$), on a 
$ [\LL: \E] > 1$. 
D'autre part, comme l'extension $\LL / \E$ est galoisienne,
\begin{equation}\label{H_p2}
 [\LL: \E] = | \Gal ( \LL / \E ) | = \text{ordre} (\tau) \leq |H_p|
\leq p.
\end{equation}
On a donc l'encadrement~:
\begin{equation} \label{encadrement} 
2 \leq [ \Q(W) : \E ([p]W) ] \leq p. 
\end{equation}  

Nous allons montrer que ce degré est exactement $p$. Pour cela
considérons une racine primitive $p^{\text{ième}}$ de l'unité, $\z_p$,
et les extensions cycliques $\Q(W,\z_p) / \Q(W)$,
$\E([p]W,\z_p) / \E([p]W) $ et $\Q([p]W,\z_p) / \Q([p]W) $ dont le
degré divise $p-1$. Remarquons que l'on a ({\it cf} (\ref{égalité})):
$$
\Q([p]W,\z_p) \subseteq \E([p]W,\z_p) \subseteq \LL([p]W,\z_p ) 
= \LL(W, \z_p) = \Q(W,\z_p). 
$$
Par ailleurs, par le lemme \ref{zeta}, l'extension $\Q(W,\z_p) /
\Q([p]W,\z_p) $ est abélienne de degré une puissance de $p$. Il en est
donc de même pour l'extension intermédiaire
$\Q(W,\z_p) / \E([p]W,\z_p) $.
Supposons que $\Q(W,\z_p) = \E([p]W,\z_p) $. On a 
$$ \E([p]W,\z_p) \subseteq \E(W,\z_p) \subseteq \LL(W, \z_p) = \Q(W,\z_p)$$
donc $\E(W,\z_p) =\E([p]W,\z_p)  $. Le lemme (\ref{zeta}) implique alors
l'existence d'un $\zz \in \ker[p] $ tel que $\E (\zz W ) = \E ([p] W )
$. Ainsi :
$$ \Q(\zz W ) \subseteq \E (\zz W ) =  \E ([p] W ) \subseteq \LL([p]W)
= \LL (W) = \Q(W).$$
L'hypothèse (\ref{6}) faite sur $W$ implique que ces quatre corps sont
égaux; en particulier $\Q(W) = \E([p] W )$, ce qui est impossible
d'après (\ref{encadrement}). 
Donc $[\Q(W,\z_p):\E([p]W,\z_p) ] = p^\alpha $ avec $\alpha \geq 1$.
Ainsi $p$ divise :
$$
  [\Q(W,\z_p) : \E([p]W) ] = [\Q(W,\z_p):\Q(W)] [\Q(W)
:\E([p]W) ] .
$$
Or $[\Q(W,\z_p):\Q(W)]$ divise $p-1$ ; le lemme de Gauss implique donc
que $p$ divise $ [\Q(W) : \E([p]W) ]$. L'encadrement
(\ref{encadrement}) nous indique alors que
$$
[ \Q(W) : \E ([p]W) ] = p. 
$$

On a ainsi montré que $[\LL : \E ] = [\Q(W) : \E([p] W ) ] = p $. On
déduit alors de (\ref{H_p2}) que le cardinal de $H_p$ est
$p$. Il en est donc de même pour $G_p$.
Notons $q:=\tilde{e}_p(\LL)$ et fixons une $q^{\text{ième}}$
racine primitive de l'unité  $\z_q = \z_m^{(m/q)}$. Comme $\LL(\z_q)
\subseteq \Q( \z_m) $, le groupe de Galois $G_p$ induit par
restriction un sous-groupe non trivial de $\Gal(\LL(\z_q)/\Q)$, qui
est nécessairement cyclique d'ordre $p$. Notons $\F$ le
corps fixé par ce sous-groupe et $\ro$ un générateur de
$\Gal(\LL(\z_q) / \F)$. Alors $\E \subseteq \F$ et
\begin{equation}\label{3.5}
\ro \z_q  =  \tilde \z_p \z_q
\end{equation}
où $\tilde \z_p$ est une racine primitive $p^{\text{ième}}$ de l'unité.

Nous allons montrer qu'il existe $\zz \in (\G_m^n)_{tors}$ tel que
$\F(\zz W)  \subseteq \F([p]W)$. Nous pouvons supposer que $\F([p]W) 
\subsetneq \F(W)$, sinon notre affirmation est triviale;
donc $\F([p]W) \subsetneq \LL(W,\z_q)$. De plus, par un argument
galoisien, on a $[\LL([p]W,\z_q) : \F([p]W) ]$ qui divise $[\LL(\z_q) :
\F]=p$. Or $ \LL([p]W,\z_q) =  \LL(W,\z_q)$, donc :
$$
[\LL(W,\z_q) : \F([p]W) ] =  [ \LL([p]W,\z_q) : \F([p]W) ] = p.
$$
En utilisant de nouveau un argument galoisien, on obtient que la
restriction~:
$$
r : \Gal ( \LL(W,\z_q) / \F([p]W) ) \rightarrow \Gal (\LL(\z_q) / \F) 
$$
est un isomorphisme de groupe.
Soit $\tilde \ro$ un générateur de $\Gal ( \LL(W,\z_q) / \F([p]W)
)$. Il existe alors  $\boldsymbol{\xi}=({\tilde \z_p}^{\alpha_1}, \dots,
{\tilde \z_p}^{\alpha_n} ) \in \ker [p] $  tel que :
$$ 
\tilde \ro \W =  \boldsymbol{\xi} W
$$
et, par (\ref{3.5}):
$$
\tilde \ro \z_q  =  \tilde \z_p \z_q.
$$
Si on pose $\zz = (\z_q^{- \alpha_1}, \dots,
\z_q^{-\alpha_n} )$, alors on a :
$$
\tilde \ro ( \zz W ) = \tilde \ro ( \zz ) \tilde \ro ( W) 
 = ({\tilde \z_p}^{-\alpha_1} \z_q^{-\alpha_1},\dots, {\tilde \z_p}^{ 
      -\alpha_n} \z_q^{-\alpha_n}) \boldsymbol{\xi} W 
 = \zz W 
$$
Ainsi $\zz W$ est stable sous l'action de $\Gal ( \LL(W,\z_q) / \F([p]W)
)$, {\it i.e.}  $\F(\zz W)  \subseteq \F([p]W)$.

On en déduit que :
$$
[ \F(\zz W) : \F ] \leq  [ \F([p] W) : \F ] \leq [\E ([p] W) : \E ] =
\delta. 
$$
Or, comme $\FF \subseteq \Q( \z_{m/p} ) $, on a $\tilde e (\F) < \tilde
e (\LL)$. On vient ainsi de contredire l'hypothèse (\ref{5}) faite sur
$W$, ce qui achève la démonstration du lemme.$\hfill \square$


 \section{Construction de la fonction auxiliaire}

 Soit $S$ un sous-espace vectoriel de $\QQ^{l}$ de dimension $d$. On
 définit la hauteur $h_2$ de $S$ comme le fait {\sc Schmidt} (voir
 \cite{Schmidt}, ch.1, \S 8) par :
 $$
 h_2(S) = \sum_v \frac{[\F_v:\Q_v]}{[\F:\Q]} \log || \xx_1 \wedge \cdots
 \wedge \xx_d ||_v
 $$ 
 où $\xx_1,\dots,\xx_d$ est une base de $S$ sur un corps de nombre $\F$
 quelconque sur lequel $S$ est rationnel, et $||\cdot||_v$ est la norme
 du sup si $v$ est ultramétrique, la norme euclidienne sinon.
 Pour $\xx \in \QQ^{l}$, nous noterons, par abus de notation,
 $h_2(\xx)$ la hauteur du sous-espace engendré par $\xx$.




Nous énonçons maintenant un théorème permettant de construire une
fonction auxiliaire dont le degré et la hauteur sont controlés.

 \begin{theo}
 Soit $V$ une hypersurface de $\G_m^n$ de degré $\oo$. 
 Soient $L$ et $T$ deux entiers naturels non nuls tels que $L \geq \oo
 T$.
 Alors, pour tout~$\eps$ strictement positif, il existe un polynôme
 non nul $F \in {\overline \Q} [X_1,\cdots,X_n] $ à coefficients
 entiers algébriques, de
 degré inférieur ou égal à $L$,
 identiquement nul sur $V$ à un ordre supérieur ou égal à $T$ tel que :
 $$h_2(F) \leq r \Big ( (T+n) \log (L+1) + L (\mm(V) + \eps) \Big ) +
 \frac{1}{2} \log {L+n \choose n} $$
 avec $ r = \frac{ {L+n \choose n}- {L - \oo T +n \choose n}} {{L- \oo
 T +n   \choose  n}}$ et où, par définition, la hauteur d'un polynôme
 est la hauteur de la famille de ses coefficients.
 \end{theo}

 {\it Démonstration. }
  Nous noterons encore $V$ l'adhérence de Zariski de $V$ dans
  $\mathbb{P}_n$.
  Soit $P \in  {\overline \Q} [X_0,\cdots,X_n] $, homogène, qui
  engendre l'idéal associé à $V$ (on a $\deg P = \oo$).
  L'ensemble 
   $E$  des polynômes homogènes de ${\overline
    \Q} [X_0,\cdots,X_n] $  de
  degré $L$ identiquement nuls à un ordre $\geq T$ sur $V$ est
  constitué des polynômes $G.P^T$ où  $G \in
  {\overline \Q} [X_0,\cdots,X_n] $ est homogène de degré $L- \oo T$. 
  Si on y ajoute le polynôme nul, c'est un ${\overline \Q}$-espace vectoriel
  de dimension ${L - \oo T +n \choose n}$. 

La démonstration est alors exactement celle du théorème 2.2 de
\cite{Am-Da03} où l'on a substitué $E \cup \{0\}$ à
$[\mathfrak{P}^{(T)}]_L$ et $H(\mathfrak{P}^{(T)};L)$ à
$ {L+n \choose n}- {L - \oo T +n \choose n}. $ Nous obtenons ainsi un
polynôme $F$ à coefficients algébriques qui satisfait les propriétés voulues.
Quitte à multiplier $F$ par un entier algébrique (ce qui ne
modifie en rien sa hauteur) , on peut supposer que $F$ est à
coefficients entiers algébriques, ce qui donne le résultat
souhaité.$\hfill \square$ 

\begin{corollaire}\label{aux}
 Soit $V$ une hypersurface de $\G_m^n$ de degré $\oo$, qui n'est pas
 de torsion. 
 Soient $L$ et $T$ deux entiers naturels non nuls tels que $ L \geq  2
 \oo T$.
 Alors il existe un polynôme
 non nul $F \in {\overline \Q} [X_1,\cdots,X_n] $ à coefficients
 entiers algébriques, de
 degré $ \leq L$, nul sur $V$ à un ordre $\geq T$ tel que :
 $$h_2(F) \leq \frac{2^{n+1} \oo T}{L} \Big ( (T+n) \log (L+1) + 2 L
  \mm(V)  \Big ) + 
 \frac{n}{2} \log (L+1) . $$
\end{corollaire}

{\it Démonstration. }
Comme $V$ n'est pas de torsion, il suffit d'appliquer le théorème
précédent avec $\eps = \mm(V)$, d'utiliser l'inégalité  ${L+n \choose
  n} \leq (L+1)^n$ et de remarquer que si $L \geq 2 \oo 
T$, on a:
$$
\frac{ {L+n \choose n}- {L - \oo T +n \choose n}} {{L- \oo
 T +n   \choose  n}} \leq \frac{2^{n+1} \oo T}{L} .
$$
$\hfill \square$

\section{Extrapolation}

Nous utiliserons à plusieurs reprises un lemme d'approximation qui
permet d'exhiber un <<dénominateur commun>> local :

\begin{lemme}\label{dénominateur}
Soient $\K$ un corps de nombres, $\OK$ son
anneau d'entiers et $v$ une place ultramétrique de $\K$.
Soit $\al=(\alpha_1,\dots,\alpha_n) \in \K^n$. Il existe alors $\beta
\in \OK$ tel que :
$$
 \begin{array}{rl}
 &(\beta \alpha_1, \dots , \beta \alpha_n) \in {\OK}^n  \\
 \text{et   }&|\beta|_v = \max \left \{1,|\alpha_1|_v,\dots,
   |\alpha_n|_v \right \}^{-1}.
 \end{array}
$$
\end{lemme}

{\it Démonstration. } Voir \cite{Am-Da99}, lemme 3.2. $\hfill \square$

\begin{pte}\label{extrapolation1}
Soit $p$ un nombre premier et soient $T_1$ et $L_1$ deux entiers
naturels non nuls. 
Supposons qu'il existe un polynôme non nul $F_1$ à coefficients
entiers algébriques, de degré au plus $L_1$, identiquement nul sur
$V_1$ avec multiplicité supérieure ou égale à $T_1$. Soit $v$ une
valuation sur $\QQ$ qui prolonge la valuation $p$-adique.
Alors, pour tout $\al \in W$ et pour tout $\tau$ qui prolonge
$\Phi_p$, on a : 
$$ |F_1^{\tau} ( \al^p ) |_v \leq p^{-T_1/e_p(\LL)} \max \left \{
   1,|\alpha_1|_v,\cdots, |\alpha_n|_v \right \}^{pL_1}.
$$
\end{pte}
 
{\it Démonstration. }
Si $p \OL = ( \pi_1 \cdots \pi_r )^{ e_p (\LL) } $, notons $Q$ l'idéal
$\pi_1 \cdots \pi_r$.
Soient $\alpha
\in W$ et $\tau \in \Gal(\QQ / \Q ) $ qui prolonge $\Phi_p$.
D'après le lemme \ref{dénominateur}, il existe une équation réduite $f
\in \OL[\X]$ de $V_1$ telle que  $|f^\tau|_v = 1 $.

Par le petit théorème de Fermat et le lemme \ref{Phi}, on a :
$$
 f(\X)^p \equiv f^\tau (\X^p) \mod Q \OL[\X].
$$
En utilisant de nouveau le lemme \ref{dénominateur}, il existe $\eta
\in \OO_{\Q(\al)}$ tel que :
$$
 \begin{array}{rl}
 &\eta \alpha_1, \dots , \eta \alpha_n \in \OO_{\Q(\al)} \\
 \text{et   }&|\eta|_v = \max \left \{1,|\alpha_1|_v,\dots,
   |\alpha_n|_v \right \}^{-1}.
 \end{array}
$$
On a alors :
$$
|\eta^{p\deg(f)} f^\tau(\al^p)|_v = |\eta^{p\deg(f)} f(\al)^p -
\eta^{p\deg(f)} 
f^\tau(\al^p) |_v \leq p^{-1/e_p(\LL)}.
$$
Donc
$$
|f^\tau ( \al^p) |_v \leq p^{-1/e_p(\LL)}\max \left
  \{1,|\alpha_1|_v,\dots,   |\alpha_n|_v \right \}^{p\deg(f)} .
$$
Comme $F_1$ est à coefficients entiers algébriques et $|f^\tau|_v=1$, nous
avons la factorisation $F_1 = q \cdot f^{T_1}$ avec $|q^\tau|_v \leq
1$. Ainsi :
$$
|F_1^\tau ( \al^p) |_v \leq p^{-T_1/e_p(\LL)}\max \left
  \{1,|\alpha_1|_v,\dots,   |\alpha_n|_v \right \}^{p L_1} .
$$
$\hfill \square$

\bigskip

Dans le cas où la ramification est <<grande>>, l'étape d'extrapolation
est différente. Nous devons au préalable établir un lemme technique.

\bigskip

\noindent {\bf Notation.}
Soit $n \in \N^*$ et $f$ un polynôme. Nous noterons $f_n$ le
polynôme dont les coefficients sont obtenus en élevant ceux de $f$ à
la puissance $n$.

\begin{lemme}\label{prod}
Soient $p$ un nombre premier, $\K$ un corps de nombres et $\OK$ son
anneau d'entiers. 
Soient $f \in \OK[\X]$ et $\zz \in \ker[p]$.
Alors  $\prod_{j=0}^{p-1} f(\zz^j \X ) \in \OK[\X]$ et
$$
 \prod_{j=0}^{p-1} f(\zz^j \X ) \equiv  f_p(\X^p)  \mod p\OK[\X].
$$
\end{lemme}

\bigskip

{\it Démonstration. }
Soit $\z_p$ une racine primitive $p$-ième de l'unité. Il existe $(a_1,
\dots, a_n) \in \N^n$ tels que $\zz = ( \z_p^{a_1}, \dots, \z_p^{a_n}) $.
On peut alors écrire le produit  $\prod_{j=0}^{p-1} f(\zz^j \X )$
comme un résultant: 
$$
\prod_{j=0}^{p-1} f(\zz^j \X ) = \text{Res}_Y (f(Y^{a_1} X_1, \dots ,
Y^{a_n} 
X_n ), Y^p - 1 ) .
$$
Ceci implique que ce produit est un polynôme à coefficients dans
$\OK$. De plus, modulo $p\OK[\X]$, on a
$$
\begin{disarray}{rcl}
\prod_{j=0}^{p-1} f(\zz^j \X )
 &\equiv & \text{Res}_Y (f(Y^{a_1} X_1, \dots ,Y^{a_n} X_n ), (Y - 1)^p
 )  \\
 &\equiv & f(\X)^p  \\[2mm]
 &\equiv & f_p(\X^p) \ ,
\end{disarray}
$$
ce qui nous donne la congruence annoncée.$\hfill \square$

\bigskip

\begin{lemme}\label{tech}
Soient $p \in \p$, $v$ une place de $\OL$ au dessus de $p$ et $\tau
\in Gal(\QQ/\Q)$ tel que $\tau_{|\LL} \in H_p$.
Il existe alors un polynôme $g \in \OL[\X]$ tel
que :
$$
\begin{array}{rl}
i) & g \  \text{est une équation réduite de  } V_p, \\
ii) & |g^\tau|_v = 1, \\
iii) & \text{Il existe $t \in \N^*$ tel que } g(\X^p) \equiv
  f_{(p^t)}(\X^{p^t}) \mod p\OL[\X], \\
 & \text{où } f \in \OL[\X] \text{ est une équation réduite de  } V_1. 
\end{array}
$$

\end{lemme}

{\it Démonstration. } 
Considérons l'ensemble algébrique suivant :
\begin{equation}\label{pVp}
[p]^{-1} V_p  = \bigcup_{\zz \in \ker [p]} \zz V_1 = \bigcup_{ \zz
  \in \HH } \zz V_1 \ 
\end{equation}
où $\HH$ est le groupe quotient défini par :
$$
\HH = \left( ^{ \ker [p]}/_{ \ker [p] \cap G_{V_1}} \right ) .
$$
(Remarquons que la variété $\zz V_1$ pour $\zz
  \in \HH$ est bien définie car celle-ci ne dépend pas, par
  construction de $\HH$, du représentant de $\zz$ choisi dans $\ker[p]$.)
Dans le dernier membre de l'égalité (\ref{pVp}), la réunion est constituée
d'hypersurfaces n'ayant aucune composante irréductible commune. En
effet, d'une part la réunion est faite modulo le stabilisateur de $V_1$,
ce qui nous assure que pour deux éléments distincts $\zz $
et $ \tilde \zz $ de $\HH$, on a $\zz 
V_1 \neq \tilde \zz V_1$; d'autre part, pour tout plongement $\sigma : \LL(W)
\hookrightarrow \QQ $ qui fixe $\LL$, on a $G_{\sigma(W)} \subseteq
G_{V_1}$, et donc de
même $\zz \sigma(W) \neq \tilde \zz \sigma(W)$. Enfin, si l'on suppose
l'existence de $\sigma$ et $\tilde \sigma$ deux plongements distincts de
$\LL(W)$ dans $\QQ$ tels que $\sigma_{|\LL} = \tilde \sigma_{|\LL}
= Id $ et pour lesquels on a $\zz \sigma(W) = \tilde \zz \tilde
\sigma(W)$ (avec $\zz,\ \tilde \zz  \in \ker[p])$, ceci implique, par
multiplication par $p$, que $\sigma([p]W) = \tilde \sigma ([p]W)$, ce
qui est impossible car $\LL(W) = \LL([p]W)$ ({\it cf} lemme \ref{i}).

Soit $f \in \OL[\X]$ une équation réduite de $V_1$ telle que $|f^{\tau}|_v=1$
({\it cf} proposition \ref{extrapolation1}); on pose 
$$
 h (\X) \ :=  \prod_{ \zz \in \HH } f( \zz \X )
$$
où l'on choisit pour chaque terme du produit un représentant de la classe
d'équivalence. Alors, d'après ce qui précède, on a $h$ est une
équation réduite pour $[p]^{-1} V_p$.
En effet, soient $s$ la dimension du
stabilisateur de $V_1$ et $\zz_1, \cdots ,\zz_{n-s}$
une base d'un supplémentaire de $\ker[p] \cap G_{V_1}$ dans
$\ker[p]$. Posons :
$$
 h (\X) \ =  \prod_{ (a_1,\dots,a_{n-s}) \in (\Z / p \Z)^{n-s} }
 f(\zz_1^{a_1}  \cdots \zz_{n-s}^{a_{n-s}} \X ) .
$$
Grâce au lemme précédent et par récurrence, on a : 
\begin{equation}\label{eq3}
h \in \OL[\X]  \quad \quad \text{ et } \quad \quad  h(X) \equiv
f_{p^{n-s}}(X^{p^{n-s}}) \mod p\OL[X]. 
\end{equation}
Par ailleurs, on a également :
$$
\prod_{\zz \in \ker[p]} f(\zz \X)  \in \OL [ \X^p ] \ ,
$$
et pour tout $\tilde{\zz} \in \ker[p]$, 
$$
\prod_{\zz \in \ker[p]\cap G_{V_1} }  f(\tilde{\zz} \zz \X ) =
(f(\tilde{\zz} \X))^{p^s}\ .
$$
On en déduit que $h(\X)^{p^s} \in \OL [ \X^p ]$. Comme $h$ n'est
pas un monôme ($h$ définit une sous-variété de $\G_m^n$), on a
également $h \in \OL [\X^p] $. On définit alors $g$ par 
$$ g(\X^p) = h (\X). $$
Ainsi, $g$ est à coefficients dans $\OL$ et est une équation réduite
de $V_p$; par construction, ayant choisi $f$ tel que
$|f^{\tau}|_v=1$, on a également $|g^{\tau}|_v=1$; enfin, comme $s
\leq n-2$, la congruence (\ref{eq3}) nous assure le point $(iii)$.
$\hfill \square$


\begin{pte}\label{extrapolation2}
Soit $p \in \p$ et soient $T_2$ et $L_2$ deux entiers
naturels non nuls. 
Supposons qu'il existe un polynôme non nul $F_2$ à coefficients
entiers algébriques, de degré au plus $L_2$, identiquement nul sur $V_p$
avec multiplicité supérieure ou égale à $T_2$. Soit $v$ une valuation
de $\QQ$ qui prolonge la valuation $p$-adique.
Alors, pour tout $\al \in W$ et pour tout $\tau$ qui prolonge un
élément $H_p$, on a : 
$$ |F_2^{\tau} ( \al^p ) |_v \leq p^{-T_2} \max \left \{
   1,|\alpha_1|_v,\cdots, |\alpha_n|_v \right \}^{pL_2}.
$$
\end{pte}

\bigskip

{\it Démonstration. }
Soient $\al \in W$ et $\tau \in \Gal(\QQ / \Q ) $ tel que $\tau_{|\LL}
\in H_p$. On suit alors exactement le même raisonnement que
dans la preuve de la proposition \ref{extrapolation1}, en substituant 
le polynôme $g$ construit dans le lemme précédent à $f$, l'idéal $(p)$
à $Q$, et en remarquant que :
$$
\begin{disarray}{rclcr}
g^\tau (\X^p) & \equiv & \left (f_{p^t} \right)^{\tau}  (\X^{p^t}) &
\mod p\OL[\X] &
\text{(lemme \ref{tech})} \\
& \equiv & f_{p^t}  (\X^{p^t}) &\mod p\OL[\X] & \text{(lemme
  \ref{H_p})} \\ 
& \equiv & f(\X)^{p^t} & \mod p\OL[\X]. &  \\ 
\end{disarray}
$$
$\hfill \square$

\section{Démonstration du théorème}

La démonstration suit le schéma classique d'une preuve de transcendance
: nous construisons tout d'abord une fonction auxiliaire s'annulant avec
forte multiplicité sur une certaine variété (étape
d'interpolation) ; puis nous en déduisons l'annulation de cette fonction sur
une variété de <<grand>> degré (étape d'extrapolation) pour aboutir à
une contradiction. 

Dans toute la suite, nous noterons $c_i$, $i \in \N$, des constantes
strictement positives ne dépendant (éventuellement) que de $n$. Nous
noterons également $C$ une constante (ne dépendant éventuellement que de $n$)
strictement positive suffisamment grande afin que les inégalités
ci-dessous soient vraies. Enfin, pour alléger les notations, on pose :
$$
 \D = \frac{\log (2 \w_{\LL} (W))  }{\log \log (5  \w_{\LL} (W)) }.
$$

Nous fixons maintenant deux paramètres :
$$
\begin{disarray}{lrcl}
&N &=& C^9  \D^5 \log ( 2 \w_{\LL} (W))  \\[1mm]
\text{et  } \ \ &E &=& C^3 \D^2 .
\end{disarray}
$$

Soit $\Lambda$ l'ensemble de nombres premiers $p \in \p $ tels que
$N/2 \leq p \leq N$. Le lemme suivant nous renseigne sur le cardinal
de $\Lambda$.

\begin{lemme}
On a 
$$ | \Lambda | \geq c_1 \frac{N}{\log C \cdot \log \log ( 5    \w_{\LL}
  (W) )} .$$
\end{lemme}  

{\it Démonstration. }
Par le théorème des nombres premiers, il existe $c_2>0$ tel que
l'ensemble des nombres premiers compris entre $N/2$ et $N$ soit de
cardinal supérieur à 
 $$
c_2 \frac{N}{\log N}
\geq
c_3 \frac{N}{\log C \cdot \log \log ( 5    \w_{\LL}
  (W) )}.
 $$
Par ailleurs, il y a au plus
$\log(|G_{V_1}/G^0_{V_1}|) / \log 2 $ premiers qui divisent
$|G_{V_1}/G^0_{V_1}|$ et
l'on a :
$$ |G_{V_1}/G^0_{V_1}| \leq \deg (G_{V_1}) \leq \deg (V_1) ^n =
\w_{\LL} (W)^n . $$
Nous pouvons ainsi majorer le cardinal de $E_{exc}(V_1)$ :
$$
|E_{exc}(V_1)| \leq \frac{n}{\log 2 } \log \w_{\LL} (W).
$$
Il existe donc une constante $c_1$  telle que le lemme soit vrai.
$\hfill \square$

\bigskip

Notons maintenant $\Lambda_1$ l'ensemble des premiers $p \in \Lambda$
tels que $e_p(\LL) \leq E$ et $\Lambda_2$ son complémentaire dans
$\Lambda$. Nous distinguons alors deux cas.

\subsection{Une majorité de premiers de $\Lambda$ sont <<peu>>
  ramifiés.}

Supposons en effet que
$$  
| \Lambda_1 | \geq \frac{c_1}{2} \frac{N}{\log C \cdot \log \log ( 5
  \w_{\LL} (W) )}.$$
Nous introduisons alors les deux nouveaux paramètres suivants :
$$
L_1 = \left [ C^8 \w_{\LL}(W) \D^6 \right ] \ \  \text{   et   } \ \
T_1 = \left [ C^4 \D^3 \right ]. 
$$

L'étape d'interpolation consiste en la construction d'une fonction
auxiliaire de petite hauteur s'annulant avec forte multiplicité sur
$V_1$ :
\begin{pte}
  Il existe un polynôme non nul $F_1$ à coefficients entiers
  rationnels de degré
  $\leq L_1$, nul à un ordre $\geq T_1$ sur $V_1$ et tel que :
  \begin{equation}\label{h_2(F_1)}
  h_2(F_1) \leq  c_4 \log C \cdot \log ( 2 \w_{\LL}(W) ) .
  \end{equation}
\end{pte}

\bigskip

{\it Démonstration. }
Remarquons tout d'abord que l'on a :
$$ \mm(V_1) = \mm(W) \ \text{ et } \deg(V_1) = \w_{\LL}(W). $$
Ainsi $L_1 \geq 2 \deg(V_1) T_1$.
Le corollaire \ref{aux} nous indique alors l'existence d'un polynôme
non nul 
$F_1$ à coefficients entiers algébriques, de degré $\leq L_1$,
s'annulant sur $V_1$ avec multiplicité $\geq T_1$, tel que :
$$
h_2(F_1)  \leq  \frac{2^{n+1} \w_{\LL}(W)  T_1}{L_1} \Big ( (T_1+n) \log
(L_1+1) + 2 L_1  \mm(W)  \Big ) +  
  \frac{n}{2} \log (L_1 +1) .
$$
Grâce aux choix des paramètres et à la majoration (\ref{1}) de
$\mm(W)$, on a :
$$
\begin{disarray}{rcl}
h_2(F_1) & \leq & c_5 \log C \cdot \log ( 2 \w_{\LL}(W) ) + c_6 C^4
\w_{\LL}(W) \D^3 \mm(W) \\[5mm]
&  \leq & c_5 \log C \cdot \log ( 2 \w_{\LL}(W) ) + c_6 C^4 c(n)
\D^{2-6(n+1)} .
\end{disarray}
$$
Si $c(n) \leq C^{-4} \log C $, on a alors :
$$
h_2(F_1) \leq  c_4 \log C \cdot \log ( 2 \w_{\LL}(W) ) .
$$
$\hfill \square$

\bigskip

La fonction auxiliaire ainsi construite s'annule alors sur des
conjugués de multiples de $W$: 
\begin{pte}
Avec les notations précédentes,
$F_1$ est nul sur les $\tau^{-1} [p] W $ pour tout $p \in
\Lambda_1$ et tout 
$\tau \in \Gal(\QQ / \Q )$ qui prolonge $\Phi_p$. 
\end{pte}

\bigskip

{\it Démonstration. }
Supposons qu'il existe un premier $p \in \Lambda_1$ et un
automorphisme $\tau \in \Gal(\QQ / \Q )$ qui prolonge $\Phi_p$ tel que
$F_1^\tau$ ne soit pas identiquement nul sur $[p]W$. 
Alors il existe $\al \in W$ de hauteur inférieure à $ 2 \mm (W) $
tel que $F_1^{\tau} (\al^p ) \neq 0$. 
Soit $\F$ un corps contenant les coefficients de $F_1^\tau$ et
$\al$. Soit $v$ une place de $\F$. Alors, par la proposition
\ref{extrapolation1},
$$
\text{si } v \mid p , \quad  
   |F_1^{\tau} (\al^p ) |_v  \ \leq \ p^{-T_1 / e_p(\LL) } \max
   \{1,|\alpha_1|_v , \dots, |\alpha_n|_v \}^{pL_1}.
$$
Par ailleurs, on a les majorations usuelles :
$$
\begin{array}{lrcl}
\text{si } v \mid \infty , &  
   |F_1^{\tau} (\al^p ) |_v &\leq &|F_1^\tau|_v (L_1+1)^n  \max
   \{1,|\alpha_1|_v , \dots, |\alpha_n|_v \}^{pL_1} \\[2mm]
\text{si } v  \nmid \infty , &  
   |F_1^{\tau} (\al^p ) |_v &\leq &\max
   \{1,|\alpha_1|_v , \dots, |\alpha_n|_v \}^{pL_1}. 
\end{array}
$$
La formule du produit donne alors :
$$
\begin{disarray}{rcl}
0 & = &\frac{1}{[\F:\Q]} \sum_{v \in \M_\F} [\F_v:\Q_v] \log | F_1^\tau
(\al^p ) |_v\\[5mm]
& \leq & -\frac{T_1}{e_p(\LL)} \log p + p L_1 h(\al) + h(F_1) + n \log (L+1).  
\end{disarray} 
$$
En utilisant le fait que $e_p(\LL) \leq E$ puis les inégalités
(\ref{h_2(F_1)}) (en remarquant que $h(F_1) \leq h_2(F_1)$) et (\ref{1}), on
obtient~:
$$
\begin{disarray}{rcl}
0& \leq & -\frac{T_1}{E} \log\frac{N}{2} 
          +2 N L_1 \mm(W)   + h_2(F_1) + n \log (L+1)\\[5mm]
& \leq & - c_7 C \log( 2 \w_{\LL}(W) ) + 2 c_8 C^{17} \w_{\LL}(W) 
          \D^{11} \log ( 2 \w_{\LL}(W) ) \mm(W) \\[2mm]
& & \hspace{30mm}  +  c_6 \log C \cdot \log ( 2 \w_{\LL}(W) ) +   c_9
\log C \cdot \log ( 2 \w_{\LL}(W) ) 
\\[5mm]
& \leq & -c_{10} C \log( 2 \w_{\LL}(W) ) + c_{8} C^{17} c(n) \D^{11-6(n+1)}
          \log ( 2 \w_{\LL}(W) ) 
\end{disarray}
$$
Si $c(n) < c_{8}^{-1} c_{10} C^{-16}$, on aboutit alors à une contradiction.
$\hfill \square$

\bigskip

On déduit de cette proposition que $F_1$ s'annule sur :
$$
\tilde {V_1} : =  \bigcup_{ p \in \Lambda_1 }
            \bigcup_{\tau :\LL([p]W) \hookrightarrow \QQ
                     \atop \tau_{|\LL} = \Phi_p}
             \tau^{-1}  [p] W   .
$$
Par ailleurs, d'après le lemme \ref{i}, si $p \in \Lambda_1$, on a :
$$[\LL([p]W) : \LL ] = [\LL(W) : \LL ] =\dd.$$
Il existe donc exactement $\dd$ morphismes distincts $\tau :
\LL([p] W ) \hookrightarrow \QQ $ qui prolongent $\Phi_p$.
De plus,
le lemme 2.3 de \cite{Am-Da99} assure que si $p \neq q$ et si $\tau_1
$ et 
$\tau_2$ sont deux éléments de $\Gal(\QQ / \Q )$, alors $\tau_1([p]W)
\neq \tau_2([q]W)$. D'où : 
$$
\begin{disarray}{rcl}
\deg \tilde{V_1} & 
= & \sum_{p \in \Lambda_1}  \dd  \deg ([p]W) \\
& = &  \sum_{p \in \Lambda_1} p^{n-1-s}\dd \deg W  \\
& \geq & |\Lambda_1 | \w_{\LL}(W)  \left (\frac{N}{2} \right
)^{n-1-s}\\[4mm]
& \geq &  c_{11}  \frac{N^{n-s}
  }{\log C \cdot \log \log ( 5  \w_{\LL} (W) )}\w_{\LL}(W)\\[4mm]
& \geq & c_{11} 
 \frac{N
  }{\log C \cdot \log \log ( 5  \w_{\LL} (W) )}\w_{\LL}(W)\\[4mm]
& \geq & c_{11}  \frac{C}{\log C} L_1,
\end{disarray}
$$ 
ce qui constitue une contradiction (pour $C$ assez grand) puisque
$F_1$ est nul sur $\tilde{V_1}$ et $\deg F_1 \leq L_1$.

\subsection{Une majorité de premiers de $\Lambda$ sont <<beaucoup>>
  ramifiés.}

Supposons maintenant que
$$  
| \Lambda_2 | \geq \frac{c_1}{2} \frac{N}{\log C \cdot \log \log ( 5
  \w_{\LL} (W) )}.$$
Nous introduisons alors les deux nouveaux paramètres suivants :
$$
L_2 = \left [ C^{9(n-s)+2} \w_{\LL}(W) \D^{6(n-s)+2} \right ] \ \
  \text{   et   } \ \
T_2 = \left [ C \D \right ]
$$
et nous considérons l'ensemble algébrique $U$, réunion d'hypersufaces :
$$
   U  =  \bigcup_{p \in \Lambda_2} V_p \ .
$$
Grâce au lemme \ref{i} :
$$
\deg U \leq |\Lambda_2| N^{n-1-s} \w_{\LL} (W) \leq c_{12} \frac{N^{n-s}
}{\log C \cdot \log \log 5 \w_{\LL} (W)}  \w_{\LL} (W)
$$
et  
$$\mm(U) =  \max_{p \in \Lambda_2} \mm ([p] \overline{W}^\LL )
  \leq N \mm(W) .$$
Comme dans le cas précédent,
\begin{pte}
  Il existe un polynôme non nul $F_2$ à coefficients entiers
  rationnels de degré 
  $\leq L_2$, nul à un ordre $\geq T_2$ sur $U$ et tel que :
  \begin{equation}\label{h_2(F_2)}
  h_2(F_2) \leq  c_{13} \log C \cdot \log ( 2 \w_{\LL}(W) ) .
  \end{equation}
\end{pte}

{\it Démonstration. }
On a $L_2 \geq 2 T_2 \cdot \deg U$. D'après le corollaire \ref{aux},
il existe donc un polynôme non nul
$F_2$ à coefficients entiers algébriques, de degré $\leq L_2$,
s'annulant sur $U$ avec multiplicité $\geq T_2$, tel que :
$$
h_2(F_2)  \leq  \frac{2^{n+1} T_2 \deg U }{L_2} \Big ( (T_2+n) \log
(L_2+1) + 2  L_2 N  \mm(W)  \Big ) +  
  \frac{n}{2} \log (L_2 +1).
$$
Grâce aux choix des paramètres et à la majoration (\ref{1}) de
$\mm(W)$, on a :
$$
\begin{disarray}{rcl}
h_2(F_2) & \leq & c_{14} \log C \cdot \log ( 2 \w_{\LL}(W) ) \\
&& \hspace{15mm}  + c_{15}
C^{1+9(n+1-s)}
\w_{\LL}(W) \D^{1+6(n-s+1)} \mm(W) \\[4mm]
&  \leq & c_{14} \log C \cdot \log ( 2 \w_{\LL}(W) ) + c_{15}
C^{1+9(n+s)}  c(n).
\end{disarray}
$$
Si $c(n) \leq C^{-(1+9(n-s+1))} \log C $, on a alors :
$$
h_2(F_2) \leq  c_{13} \log C \cdot \log ( 2 \w_{\LL}(W) ) .
$$
$\hfill \square$

\bigskip

Passons maintenant à l'extrapolation :
\begin{pte}
Avec les notations précédentes,
$F_2$ est nul sur les $\tau [p] W $ pour tout $p \in
\Lambda_2$ et tout $\tau \in \Gal(\QQ / \Q )$ tel que $\tau_{|\LL} \in
H_p$. 
\end{pte}

{\it Démonstration. }
Supposons qu'il existe un premier $p \in \Lambda_2$ et un
automorphisme $\tau \in \Gal(\QQ / \Q )$ qui prolonge un élément de
$H_p$ tel que $F_2^\tau$ ne soit pas identiquement nul sur $[p]W$. 
Alors il existe $\al \in W$ de hauteur inférieure à $2 \mm (W)$
tel que $F_2^{\tau} (\al^p ) \neq 0$. 
Soit $\F$ un corps contenant les coefficients de $F_2^\tau$ et
$\al$. Soit $v$ une place de $\F$. Alors, par la proposition
\ref{extrapolation2},
$$
\text{si } v \mid p , \quad  
   |F_2^{\tau} (\al^p ) |_v  \ \leq \ p^{-T_2 } \max
   \{1,|\alpha_2|_v , \dots, |\alpha_n|_v \}^{pL_2}.
$$
La formule du produit donne alors :
$$
\begin{disarray}{rcl}
0 & = &\frac{1}{[\F:\Q]} \sum_{v \in \M_\F} [\F_v:\Q_v] \log | F_2^\tau
(\al^p ) |_v\\[5mm]
& \leq & -T_2 \log p + p L_2 h(\al) + h(F_2) + n \log (L+1).  
\end{disarray} 
$$
En utilisant les inégalités
(\ref{h_2(F_2)}) et (\ref{1}), on
obtient~:
$$
\begin{disarray}{rcl}
0& \leq & -T_2 \log\frac{N}{2} 
          + 2 N L_2 \mm(W)  + h_2(F_2) + n \log (L+1)\\[2mm]
& \leq & - c_{17} C \log( 2 \w_{\LL}(W) ) \\
&& \hspace{5mm}+ c_{18}\w_{\LL}(W) C^{2+9(n-s+1)} 
          \D^{1+6(n-s+1)} \log \log (5  \w_{\LL}(W) ) 
          \mm(W) \\[2mm]
& \leq & - c_{17} C \log( 2 \w_{\LL}(W) ) + c_{18} C^{2+9(n+1)} c(n) 
\log \log (5  \w_{\LL}(W) ) .
$$
\end{disarray}
$$
Ce qui constitue une contradiction dès que $c(n) < {c_{18}}^{-1} c_{17}
C^{-(1+9(n+1))}$.$\hfill \square$

\bigskip

On déduit de cette proposition que $F_2$ s'annule sur :
$$
\tilde {U} : =  \bigcup_{ p \in \Lambda_2 }
            \bigcup_{\tau :\LL([p]W) \hookrightarrow \QQ
                     \atop \tau_{|\LL} \in H_p}
             \tau  [p] W   .
$$

Soit $p \in \Lambda_2$ et soit $\sigma \in H_p$.
Comme dans le cas précédent, il existe exactement $\dd$ morphismes
distincts $\tau : \LL([p] W ) \hookrightarrow \QQ $ qui prolongent
$\sigma$. 
De plus, la dernière assertion du lemme \ref{H_p} assure que si
$\tilde{\sigma} \in H_p \backslash \{\sigma \} $ et si $\tau$ et $\tilde \tau
\in \Gal(\QQ / \Q )$ sont tels que $\tau_{|\LL} = \sigma$
et ${\tilde \tau}_{|\LL} = \tilde \sigma$ , alors $\tau ([p]W)
\neq \tilde {\tau}([p]W)$. Associant cela au fait que si $p \neq q$ et
si $\tau_1,\ \tau_2 \in \Gal(\QQ / \Q )$, alors $\tau_1([p]W) \neq
\tau_2([q]W)$, on obtient : 
$$
\begin{disarray}{rcl}
\deg \tilde{U} & 
= & \sum_{p \in \Lambda_2}  \dd |H_p|  \deg ([p]W) \\
& = & \w_{\LL}(W) \sum_{p \in \Lambda_2} |H_p| p^{n-1-s} .
\end{disarray}
$$
Remarquons enfin que :
$$
|H_p| \geq \max ( p, e_p(\LL) ) \geq \max (N/2 ,E) \geq E.
$$
On a alors :
$$
\begin{disarray}{rcl}
\deg \tilde{U}   & \geq & \w_{\LL}(W) |\Lambda_2 | E  \left
  (\frac{N}{2} \right )^{n-1-s}  \\[4mm]
& \geq &  c_{19}  \frac{N^{n-s}
  }{\log C \cdot \log \log ( 5  \w_{\LL} (W) )} E \w_{\LL}(W)\\[4mm]
& \geq & c_{20} \frac{C}{\log C} L_2,
\end{disarray}
$$ 
ce qui est de nouveau une contradiction et achève la preuve du
théorème.

\bibliographystyle{alpha}
\bibliography{Biblio.bib}

\begin{thebibliography}{Zha95b}

\bibitem[AD99]{Am-Da99}
F.~Amoroso and S.~David.
\newblock {<< Le probl\`eme de Lehmer en dimension sup\'erieure>>}.
\newblock {\em J. Reine Angew. Math.}, 513:145--179, 1999.

\bibitem[AD00]{Am-Da00}
F.~Amoroso and S.~David.
\newblock {<< Minoration de la hauteur normalisée des hypersurfaces>>}.
\newblock {\em Acta Arith.}, 92(4):339--366, 2000.

\bibitem[AD03]{Am-Da03}
F.~Amoroso and S.~David.
\newblock {<< Minoration de la hauteur normalis\'ee dans un tore>>}.
\newblock {\em J. Inst. Math. Jussieu}, 2(3):335--381, 2003.

\bibitem[AD04]{Am-Da04}
F.~Amoroso and S.~David.
\newblock {<< Distribution des points de petite hauteur dans les groupes
  multiplicatifs>>}.
\newblock {\em Ann. Sc. Norm. Sup. Pisa, Cl. Sci (5)}, III(2):325--348, 2004.

\bibitem[AZ00]{Am-Za}
F.~Amoroso and U.~Zannier.
\newblock {<< A relative Dobrowolski lower bound over abelian extensions>>}.
\newblock {\em Ann. Sc. Norm. Sup. Pisa, Cl. Sci (4)}, XXIX(3):711--727, 2000.

\bibitem[Hin88]{Hindry}
M.~Hindry.
\newblock {<<Autour d'une conjecture de S. Lang>>}.
\newblock {\em Invent. Math.}, 94:575--603, 1988.

\bibitem[Phi91]{Ph1}
P.~Philippon.
\newblock {<< Sur des hauteurs alternatives I>>}.
\newblock {\em Math. Ann.}, 289(2):255--283, 1991.

\bibitem[Phi94]{Ph2}
P.~Philippon.
\newblock {<< Sur des hauteurs alternatives II>>}.
\newblock {\em Ann. Inst. Fourier}, 44(4):1043--1065, 1994.

\bibitem[Phi95]{Ph3}
P.~Philippon.
\newblock {<< Sur des hauteurs alternatives III>>}.
\newblock {\em J. Math. Pures Appl.(9)}, 74(4):345--365, 1995.

\bibitem[Sch91]{Schmidt}
W.~M. Schmidt.
\newblock {\em {<<Diophantine approximations an diophantine equations>>}}.
\newblock Lecture Notes in Mathematics, 1467, Springer-Verlag, 1991.

\bibitem[Zha92]{zhang1992}
S.~Zhang.
\newblock {<< Positive line bundles on arithmetic surfaces>>}.
\newblock {\em Ann. Math.}, 136(3):569--587, 1992.

\bibitem[Zha95a]{zhang1995a}
S.~Zhang.
\newblock {<< Positive line bundles on arithmetic varieties >>}.
\newblock {\em J. Am. Math. Soc.}, 8(1):187--221, 1995.

\bibitem[Zha95b]{zhang1995b}
S.~Zhang.
\newblock {<< Small points and adelic metrics>>}.
\newblock {\em J. Algebr. Geom.}, 4(2):281--300, 1995.

\end{thebibliography}

\end{document}